\documentclass[12pt]{article}
\usepackage{fullpage}
\usepackage{amsmath, amsthm, amssymb}
\numberwithin{equation}{section}
\begin{document}
\author{Ajai Choudhry and Arman Shamsi Zargar}
\title{The diophantine equation $x^4+y^4=z^4+w^4$}
\date{}
\maketitle

\begin{abstract}
Since 1772, when Euler first described two methods of obtaining  two pairs of biquadrates with equal sums, several methods of solving the diophantine equation $x^4+y^4=z^4+w^4$ have been published. All these methods yield parametric solutions  in terms of homogeneous bivariate polynomials of odd degrees. In this paper we describe a method that yields three parametric solutions of the aforesaid diophantine equation in terms of homogeneous bivariate polynomials of even  degrees, namely degrees~$74$, $88$ and $132$ respectively.
\end{abstract}



\section{Introduction}
Ever since 1772, when Euler found two pairs of biquadrates with equal sums, the diophantine equation
\begin{equation} \label{quartic}
	x^4+y^4=z^4+w^4,
\end{equation}
has been the subject of extensive investigations by a number of mathematicians (see for instance, \cite{Bernstein, Brudno, Choudhry}, \cite[pp.~644--647]{Dickson}, \cite{Dyer, Lander2, Lander3, Steggall, S-Dyer, Wroblewski, Zajta}). Euler solved Eq.~\eqref{quartic} by two methods that lead to parametric solutions giving the values of $x, y, z$ and $w$ in terms of homogeneous bivariate polynomials of degrees~$7$ and $13$ (see~\cite[p.~1062]{Lander1}), respectively. Subsequently, several other parametric solutions of degrees~$11$, $13$, $19$, $21$, $31$ and higher degrees have been obtained. It is noteworthy that all the known parametric solutions of Eq.~\eqref{quartic} are of odd degree.

Guy \cite[pp.~212--213]{Guy} refers to  two methods of generating new solutions  starting from  known solutions of Eq.~\eqref{quartic}.  He writes that these two methods ``generate all nonsingular parametric solutions, i.e., all solutions which correspond to points on curves with no singular points", and adds that all nonsingular solutions have odd degree. He further states that singular solutions do exist, and some of the singular solutions have even degree . While two and a half centuries have passed since Euler found the first solution of Eq.~\eqref{quartic}, and parametric solutions of even degree are known to exist, no such solution has been found till now. 
  
In this paper we obtain three parametric solutions of Eq.~\eqref{quartic} of even degrees, namely degrees~$74$, $88$ and $132$. We first describe, in Section~\ref{newsols}, a method of generating two new solutions starting from a known solution of Eq.~\eqref{quartic} and in Section~\ref{section3}, we apply the method to obtain parametric solutions of even degrees. We conclude the paper by mentioning some open problems related to the quartic~Eq.~\eqref{quartic}.

All computations for this paper were  carried out using the software MAPLE. 

\section{Generating new parametric solutions from old}\label{newsols}
To generate new solutions of Eq.~\eqref{quartic} starting from a known one, we will follow closely both the method and the notation of Richmond~\cite{Richmond} who has described how new solutions of the diophantine equation 
\begin{equation} \label{Rquartic}
	ax^4+by^4+cz^4+dw^4=0, 
\end{equation}
may be obtained from a known solution if  $a, b, c, d$ are integers such that their product $abcd$ is a perfect square. In the case of Eq.~\eqref{quartic}, we have $(a, b, c, d)=(1, 1, -1, -1)$,  hence the product $abcd$ is a square and we are assured of a rational solution on applying Richmond's method.

Richmond has adopted a geometrical approach to solve the problem. Let $P$ be a known rational point on the quartic surface defined by Eq.~\eqref{Rquartic}. The tangent plane at $P$ intersects the surface~\eqref{Rquartic} in a quartic curve having a double point at $P$,  and any line in the tangent plane passing through $P$ intersects the surface~\eqref{Rquartic} in two additional  points. When such a line is  a tangent to the curve at $P$, that is, an inflectional tangent of the surface, one of the additional points coincides with the point $P$ and we can find the fourth point of intersection if the equation of the tangent can be determined.

If the known rational point $P$ on the surface~\eqref{quartic} has coordinates $(X, Y, Z, W)$, then
\begin{equation}
X^4+Y^4-Z^4-W^4 =0. \label{condP}
\end{equation} 

Following Richmond, we will denote the coordinates of any other point $Q$ as $(pX, qY,$ $ rZ, sW)$. If the point $Q$ lies on the tangent plane of $P$, then
\begin{equation}
X^4p+Y^4q-Z^4r-W^4s=0. \label{condQ}
\end{equation}
Moreover, if $Q$ lies on an inflectional tangent of $P$, then
\begin{equation}
X^4p^2+Y^4q^2-Z^4r^2-W^4s^2=0. \label{condQinflex}
\end{equation}

If $(p, q, r, s)$ is any solution of the simultaneous equations~\eqref{condQ} and \eqref{condQinflex}, then adding a constant to each of $p, q, r$ and $s$, yields another solution. Hence, we may impose the  condition,
\begin{equation}
p+q+r+s=0. \label{condpqrs}
\end{equation}
Richmond has shown that the three equations~\eqref{condQ}, \eqref{condQinflex} and \eqref{condpqrs} together imply that
\[
X^4p^4+Y^4q^4-Z^4r^4-W^4s^4=0,
\]
and hence we obtain a new point $(pX, qY, rZ, sW)$ on the quartic surface~\eqref{quartic}. 

We will now solve the simultaneous~Eqs.~\eqref{condQ}, \eqref{condQinflex} and \eqref{condpqrs} to obtain two solutions for $(p, q, r, s)$ and thus obtain two new solutions of Eq.~\eqref{quartic} starting from the known solution $(X, Y, Z, W)$.

On solving the two linear equations~\eqref{condQ} and \eqref{condpqrs} for $r$ and $s$, we get,
\begin{equation}
\begin{aligned}
r & = -\big((X^4 + W^4)p + (W^4 + Y^4)q\big)/(W^4 - Z^4),\\
s & = \big((X^4 + Z^4)p + (Y^4 + Z^4)q\big)/(W^4 - Z^4).
\end{aligned}
\label{valrs}
\end{equation}
On substituting these values of $r$ and $s$ in Eq.~\eqref{condQinflex}, we get
\begin{multline}
\{(W^4 + Z^4)X^8 - (W^8 - 6W^4Z^4 + Z^8)X^4 + W^4Z^4(W^4 + Z^4)\}p^2 \\
+ \{(2W^4Y^4 + 4W^4Z^4 + 2Y^4Z^4)X^4 + 2W^4Z^4(W^4 + 2Y^4 + Z^4)\}pq \\
+ \{(W^4 + Z^4)Y^8 - (W^8 - 6W^4Z^4 + Z^8)Y^4 + W^4Z^4(W^4 + Z^4)\}q^2=0. \label{condQinflex2}
\end{multline}

Eq.~\eqref{condQinflex2} is a quadratic equation in $p$, and its discriminant with respect to $p$, using the relation~\eqref{condP}, works out to $64X^4Y^4Z^4W^4(W^4 - Z^4)^2q^2$. Accordingly, we get the following two rational solutions of Eq.~\eqref{condQinflex2}:
\begin{multline}
p = -q\{(W^4Y^4 + 2W^4Z^4 + Y^4Z^4)X^4+4Y^2Z^2W^2(W^4-Z^4)X^2 \\
+Z^4W^4(W^4+2Y^4+Z^4)\}/\{(W^4 + Z^4)X^8 \\ 
- (W^8 - 6W^4Z^4 + Z^8)X^4 + Z^4W^4(Z^4+W^4)\},
\nonumber
\end{multline}
and
\begin{multline}
p = -q\{(W^4Y^4 + 2W^4Z^4 + Y^4Z^4)X^4-4Y^2Z^2W^2(W^4-Z^4)X^2\\
+Z^4W^4(W^4+2Y^4+Z^4)\}/\{(W^4 + Z^4)X^8 \\
- (W^8 - 6W^4Z^4 + Z^8)X^4 + Z^4W^4(Z^4+W^4)\}.
\nonumber
\end{multline}

Now, using the relations~\eqref{valrs}, we get two sets of values of $p, q, r, s$ satisfying the simultaneous~Eqs.~\eqref{condQ}, \eqref{condQinflex} and \eqref{condpqrs} leading to two new solutions $(X_1, Y_1, Z_1, W_1)$ and $(X_2, Y_2, Z_2, W_2)$ of the diophantine 
equation~\eqref{quartic}. These new solutions of Eq.~\eqref{quartic} may be written explicitly as follows:
\begin{equation}
\begin{aligned}
X_1 & = X\big((W^4Y^4 + 2W^4Z^4 + Y^4Z^4)X^4 + 4W^2X^2Y^2Z^2(W^4-Z^4)\\
 & \quad \quad + W^4Z^4(W^4 + 2Y^4 + Z^4)\big),\\
Y_1 & = Y\big((W^4 + Z^4)X^8 - (W^8 - 6W^4Z^4 + Z^8)X^4 + W^4Z^4(W^4 + Z^4)\big),\\
Z_1 & = Z\big(W^4X^8 - 4W^2Y^2Z^2X^6 - (W^8 + 2W^4Y^4 - 2W^4Z^4 - Y^4Z^4)X^4\\
  & \quad \quad- 4W^6Y^2Z^2X^2 - W^4Y^4Z^4\big),\\
W_1 &= W\big(Z^4X^8 + 4W^2Y^2Z^2X^6 + (W^4Y^4 + 2W^4Z^4 - 2Y^4Z^4 - Z^8)X^4\\
  & \quad \quad + 4W^2Y^2Z^6X^2 - W^4Y^4Z^4\big),
\label{newsol1}
\end{aligned}
\end{equation}
and
\begin{equation}
\begin{aligned}
X_2 & =X\big((W^4Y^4 + 2W^4Z^4 + Y^4Z^4)X^4 - 4W^2X^2Y^2Z^2(W^4-  Z^4)\\
 & \quad \quad  + W^4Z^4(W^4 + 2Y^4 + Z^4)\big),\\
Y_2 & =Y \big((W^4 + Z^4)X^8 - (W^8 - 6W^4Z^4 + Z^8)X^4 + W^4Z^4(W^4 + Z^4)\big),\\
Z_2 & = Z\big(W^4X^8 + 4W^2Y^2Z^2X^6 - (W^8 + 2W^4Y^4 - 2W^4Z^4 - Y^4Z^4)X^4\\
 & \quad \quad  + 4W^6Y^2Z^2X^2 - W^4Y^4Z^4\big),\\
W_2 &= W\big(X^8Z^4 - 4W^2Y^2Z^2X^6 + (W^4Y^4 + 2W^4Z^4 - 2Y^4Z^4 - Z^8)X^4\\
  & \quad \quad - 4W^2Y^2Z^6X^2 - W^4Y^4Z^4\big).
\label{newsol2}
\end{aligned}
\end{equation}

With the values of $X_i, Y_i, Z_i, W_i, i=1, 2$, being given by \eqref{newsol1} and \eqref{newsol2}, we have verified, both for $i=1$ and $i=2$, that $X_i^4+Y_i^4-Z_i^4-W_i^4$ has a factor $(X^4+Y^4-Z^4-W^4)$ confirming the fact that when $(X, Y, Z, W)$ is a known solution of the diophantine~Eq.~\eqref{quartic}, two new solutions of \eqref{quartic} are given by \eqref{newsol1} and \eqref{newsol2}.

If $(X, Y, Z, W)$ is a known solution of Eq.~\eqref{quartic}, in addition to the 
solutions~\eqref{newsol1} and \eqref{newsol2}, we may try to obtain new solutions of 
Eq.~\eqref{quartic} by taking the initial known solution as any other solution of \eqref{quartic} obtained by changing the signs of $X, Y, Z$ and $W$, or, in view of the symmetry of Eq.~\eqref{quartic}, by interchanging $X$ and $Y$ or by interchanging $W$ and $Z$. It has, however, been verified that all the solutions thus obtained are equivalent to one of the two solutions~\eqref{newsol1} and \eqref{newsol2}. Thus, applying the method described above, we can obtain only two new solutions of Eq.~\eqref{quartic} from any given solution $(X, Y, Z, W)$.

If we take the trivial solution of Eq.~\eqref{quartic} as the initial known solution, we get only trivial solutions of Eq.~\eqref{quartic}. However, starting from a nontrivial solution of \eqref{quartic}, we do obtain significant new solutions of \eqref{quartic}.

\section{Even degree parametric solutions of Eq.~\eqref{quartic}}\label{section3}
We will now obtain new solutions of Eq.~\eqref{quartic} by taking $(X, Y, Z, W)$ as a nontrivial parametric solution of \eqref{quartic}. The published parametric solutions of 
Eq.~\eqref{quartic} include one solution of degree~$7$~\cite[p.~260]{Hardy-Wright}, one of degree~$11$~\cite{Zajta}, two solutions of degree~$13$~(\cite{Brudno, Zajta}), two of 
degree~$19$~\cite{Brudno}, one solution of degree~$21$~\cite{Choudhry}, and one of degree~$31$~\cite{Brudno}. We obtained the two new 
solutions~\eqref{newsol1} and \eqref{newsol2} of Eq.~\eqref{quartic} by taking each of the published solutions of degrees~$\leq 21$ as the initial known solution. While most of the new parametric solutions thus obtained were of odd degree, we were able to obtain  three solutions of even degrees, namely degrees~$74$, $88$ and $132$. 

The solution of degree~$74$ was obtained by taking the initial solution  $(X, Y, Z, W)$ of 
Eq.~\eqref{quartic} as the parametric solution of degree~$11$ given by Zajta~\cite[p.~651]{Zajta}. For brevity, we will denote a polynomial $c_0u^n+c_1u^{n-1}v+ \cdots + c_nv^n$ by $(c_0, c_1, \ldots, c_n)$. Accordingly, we may write the initial known solution $(X, Y, Z, W)$ of degree $11$ as follows:
\begin{equation*}
\begin{aligned}
X & = (-1, -1, 4, 17, 33, 49, 58, 52, 32, 12, 2, 0), \\
Y & = (1, 4, 8, 7, 5, 17, 44, 64, 58, 34, 12, 2),\\
Z & = (1, 3, 8, 13, 9, -13, -44, -64, -58, -34, -12, -2),\\
W & = (1, 2, 2, 7, 27, 59, 78, 66, 36, 12, 2, 0).
\end{aligned}
\end{equation*}

With these values of  $(X, Y, Z, W)$, on using the relations~\eqref{newsol1} and \eqref{newsol2}, we get two new  parametric solutions of Eq.~\eqref{quartic} of degrees~$71$ and $74$ respectively. While the solution of degree~$71$ is not of much interest, this is the first time that a solution of even degree of the quartic~equation~\eqref{quartic} has been found. Accordingly, even though the solution is cumbersome to write, we give below the solution of degree~$74$:
	\scriptsize
\[
\begin{aligned}
	&&x=( 1, && 17, && 152, 
	&& 895, && 3677,\\ && 10355, 
	&& 16752, && -6222, && -140835, 
	&& -557833,\\ && -1711920, &&  -5515487, 
	&& -19959999, && -72345519, && -238033554, \\
	&& -685646800, && -1714953501, && -3721974065, 
	&& -7015296151, && -11517835142,\\ && -16609064411,
	&& -21365615196, && -24731053881, && -23868601209, 
	&& -8693923546,\\ && 47765761580, && 205456220674, 
	&& 601739916040, && 1591967993962, && 4095994601060,\\ 
	&& 10211819410484, && 23927266424288, && 51389141908972, 
	&& 99959111822592, && 175616228744844,\\ && 279322120088436, 
	&& 404232436053368, && 535981332047440, && 656753805669640, 
	&& 750917495762848,\\ && 807744919686856, && 817770961417232, 
	&& 764169978445936, && 616307330131792, && 334329090096416,\\ 
	&& -111838755728688, && -716177198849232, && -1418334840709488, 
	&& -2107891661329696, && -2654267430853920,\\ && -2950439684299520, 
	&& -2949879928244256, && -2679452870276352, && -2224041503518848, 
	&& -1692791544749664,\\ && -1183811970766368, && -761319578137344, 
	&& -450268872381568, && -244732999149312, && -122074267257536, \\
	&& -55765072370368, && -23264034935872, && -8831153717696, 
	&& -3036512628032, && -940331642368,\\ && -260402502784, 
	&& -63912753024, && -13745943296, && -2552757120, 
	&& -401389312,\\ && -52002816, && -5334272, 
	&& -406528, && -20480, && -512),
	\end{aligned}
	\]
	\[
	\begin{aligned}
	&& y= (1, && 12, && 68, 
	&& 259, && 1071,\\ && 6503,  
	&& 39416, && 188722, && 700231, 
	&& 2055472,\\ && 4846760, && 9153077, 
	&& 13239607, && 12076495, && -1384546,\\ 
	&& -22264022, && 15596517, && 361584212, 
	&& 1706206161, && 5768715671,\\ && 16854139347, 
	&& 45862872154, && 118757789579, && 289692118217, 
	&& 652845273382,\\ && 1338294128486, && 2475161995622, 
	&& 4122489145244, && 6204532963276, && 8511520885474,\\ 
	&& 10788548220988, && 12827574605276, && 14416231656980, 
	&& 15197532538952, && 15046142552324,\\ && 16116160860964, 
	&& 27435630089032, && 71120958059592, && 186419869959416, 
	&& 425938058082816, \\&& 840169049393216, && 1452600163713336, 
	&& 2235386195697664, && 3099438534892656, && 3908564730337520, \\
	&& 4515911363877424, && 4808784908497216, && 4742686437811904, 
	&& 4350587246441184, && 3725781657047328,\\ && 2988512826804672, 
	&& 2251577528136096, && 1596977048727968, && 1067985817108224, 
	&& 673860524037920,\\ && 400989088511520, && 224684399069824, 
	&& 118226429889024, && 58200428385216, && 26679801745344,\\ 
	&& 11327184587264, && 4426678166144, && 1581511834432, 
	&& 512578395840, && 149388715136,\\ && 38749162752, 
	&& 8834393088, && 1742949248, && 291587328, 
	&& 40234112,\\ && 4398336, && 357376, 
	&& 19200, && 512, && 0),
\end{aligned}
\]

\[
\begin{aligned}
	&&z=(1, && 11, && 44, 
	&& -43, && -1521,\\ && -10349, 
	&& -48094, && -184582, && -626261, 
	&& -1885097, \\&& -4862686, && -10139275, 
	&& -15292077, && -10609211, && 20915604, \\
	&& 89748564, && 185154453, && 355944759, 
	&& 1120036171, && 4531415148,\\ && 16199711507, 
	&& 48596287898, && 126583492119, && 297031705531, 
	&& 643010863626,\\ && 1294952571722, && 2416410456570, 
	&& 4133999002140, && 6405227285108, && 8887453558722,\\ 
	&& 10949047727476, && 11952001853080, && 11732733780020, 
	&& 10825910502168, && 9623681783220,\\ && 5699607762540, 
	&& -10738885867128, && -61883427411640, && -184858984506696, 
	&& -427986591320912,\\ && -838134941620800, && -1440124815325576, 
	&& -2214374565389312, && -3083234892705296, && -3916367237113744, \\
	&& -4559139881417136, && -4877129937251328, && -4800033216891168, 
	&& -4345979727260704, && -3615037502758336, \\&& -2754878085393536, 
	&& -1914161146105536, && -1203148630872352, && -674978585543232, 
	&& -329583897465504,\\ && -132419348608032, && -36557304494592, 
	&& 666392821888, && 9754035865536, && 8470935522944,\\ 
	&& 5115912504576, && 2515315350400, && 1054808007744, 
	&& 383790993088, && 121767258240,\\ && 33629791232, 
	&& 8030510080, && 1638885248, && 280917248,  
	&& 39421056,\\ && 4357376, && 356352, 
	&& 19200, && 512, && 0), 
	\end{aligned}
	\]

	\[
	\begin{aligned}
	&&w=(1, && 18, && 158, 
	&& 891, && 3557,\\ && 10197, 
	&& 19042, && 11112, && -51283, 
	&& -81800,\\ && 792900, && 6076883, 
	&& 25876793, && 85316469, && 244795852, \\
	&& 646730252, && 1593093855, && 3592501192, 
	&& 7198043589, && 12441224263,\\ && 18044751915, 
	&& 21420089180, && 20943326583, && 20712001441, 
	&& 34952512262,\\ && 89353280986, && 229583243194, 
	&& 578562542448, && 1518227607774, && 4064995267080,\\ 
	&& 10413927381044, && 24448963869396, && 51846443160436, 
	&& 99369704162752, && 173128682714556,\\ && 275899919403068,  
	&& 403984399138728, && 544562575159256, && 675160917410712, 
	&& 767365027699696,\\ && 796346654381272, && 755106814522256,  
	&& 667938723818320, && 594168862898784, && 614577345224672,\\ 
	&& 800741879613040, && 1178674777272112, && 1705240992660656, 
	&& 2272595710082912, && 2742355103542944,\\ && 2995060811603392, 
	&& 2971951766681280, && 2690549300561088, && 2229865867200576, 
	&& 1695948860452832,\\ && 1185506921484768, && 762180796885504, 
	&& 450669693171776, && 244900386259584, && 122136141687808, \\
	&& 55785094561728, && 23269642319168, && 8832494068800, 
	&& 3036781077824, && 940375537408,\\ && 260408140288, 
	&& 63913287296, && 13745976576, && 2552758144, 
	&& 401389312,\\ && 52002816, && 5334272,  
	&& 406528, 	&& 20480, && 512).
\end{aligned}
\]
\normalsize

Similarly, we obtained the parametric solution of degree~$88$ using the relations~\eqref{newsol2} and starting from the following solution of degree~$13$ given by Zajta~\cite[p.~651]{Zajta}:
\[
\begin{aligned}
X & = (1, 3, 10, 22, 44, 67, 88, 95, 84, 58, 30, 10, 2, 0),\\
Y & = (0,0,3,9,24,45,72,91,94,80,54,28,10,2),\\
Z & = (1,3,10,22,40,63,82,95,94,80,54,28,10,2),\\
W & = (0,2,5,15,28,47,64,73,66,48,26,10,2,0).
\end{aligned}
\]

Finally, we obtained the parametric solution of degree~$132$ using the relations~\eqref{newsol2} and starting from the following solution of degree~$19$ given by Brudno~\cite{Brudno}:
\[
\begin{aligned}
X & =(0, 1, 3, -15, 15, 6, -45, 82, -15, -123, 171, -159, 159, -98, 30, -12, 0, 3, 0, 1),\\
Y& = (1, -1, -3, -3, 21, -12, -44, 86, -93, 87, 3, -135, 142, -100, 72, -36, 12, -9, 1, \\
&\quad \quad -1),\\
Z & = (1, -1, -3, -3, 21, -6, -44, 62, 15, -129, 165, -129, 88, -46, 18, -6, 12, -3, 1, \\
&\quad \quad -1),\\
W & = (1, -3, 3, 21, -60, 27, 58, -75, 57, -63, 63, -87, 100, -66, 36, -18, 9, 0, 1).
\end{aligned}
\]

As the solutions of degrees~$88$ and $132$ are too cumbersome to write, we do not give them explicitly. Any interested reader can readily find these solutions using a symbolic algebra computation software, the related values of $(X, Y, Z, W)$ given above and the 
relations~\eqref{newsol2}. 

\section{Some open problems}\label{openprob}
Parametric solutions of even degree of the classical diophantine equation~\eqref{quartic} had remained frustratingly elusive till now. In this paper we obtained three parametric solutions of even degrees~$74$, $88$ and $132$ respectively. It would be useful  to find even degree solutions of lower degrees. In fact, it would be interesting to determine the parametric solution or solutions which have the least even degree.

It also remains to be investigated whether the even degree parametric solutions yield any new numerical solutions of Eq.~\eqref{quartic} that are not  generated by the solutions of odd degree. 

\bibliographystyle{amsplain}

\noindent Postal Address: Ajai Choudhry, \\
\noindent \hspace{1.05 in} 13/4 A Clay Square, \\
\noindent \hspace{1.05 in} Lucknow-226001, INDIA.

\noindent E-mail address: ajaic203@yahoo.com

\medskip

\noindent Postal Address: Arman Shamsi Zargar, \\
\noindent \hspace{1.05 in} Department of Mathematics and Applications,\\
\noindent \hspace{1.05 in} University of Mohaghegh Ardabili, Ardabil, IRAN.

\noindent E-mail address: zargar@uma.ac.ir

\end{document}